\documentclass[11pt]{article}
\usepackage{graphicx,epic,eepic}
\usepackage{amsmath, amsthm, verbatim, amsfonts, amssymb, color}
\usepackage{mathrsfs}
\usepackage{dsfont}
\usepackage[a4paper,margin=3cm]{geometry}
\usepackage{epstopdf}

\numberwithin{equation}{section}
\renewcommand\labelenumi{\textup{\alph{enumi})}}
\renewcommand\theenumi\labelenumi

\newcommand\nat{\mathds{N}}

\newcommand\eup{\mathrm{e}}

\usepackage{marginnote}

\usepackage{amssymb,amsmath,amsthm,mathdesign}
\usepackage{pdfsync}
\usepackage{color}
\usepackage[colorlinks]{hyperref}

\newtheorem{theorem}{Theorem}[section]
\newtheorem{lemma}[theorem]{Lemma}

\newtheorem{proposition}[theorem]{Proposition}

\newtheorem{assumption}[theorem]{Assumption}
\newtheorem{remark}[theorem]{Remark}

\numberwithin{equation}{section}

\newcommand{\be}{\begin{equation}}
\newcommand{\ee}{\end{equation}}

\newcommand{\bes}{\begin{equation*}}
\newcommand{\ees}{\end{equation*}}

\def\E{\bE}

\def\P{\bP} 

\newcommand{\N}{\mathbf{varrho}}

\newcommand{\R}{\mathbf{R}}

\renewcommand{\d}{{\rm d}}

\renewcommand{\geq}{\geqslant}
\renewcommand{\leq}{\leqslant}
\renewcommand{\ge}{\geqslant}
\renewcommand{\le}{\leqslant}

\renewcommand{\P}{\mathrm{P}}

\newcommand{\rd}{{\mathbb R^d}}

\def\R{{\mathbb R}}
\def\N{{\mathbb N}}

\def\a{\alpha}

\def\E{{\mathbb E}}
\def\P{{\mathbb P}}

\def\m1{\mathbf{1}}

\allowdisplaybreaks[1]

\allowdisplaybreaks[1]

\title{\textbf{Pathwise Blowup of space-time fractional SPDEs}}
\author{	\textbf{Chang-Song Deng} \\School of Mathematics and Statistics, Wuhan University, Wuhan 430072, China\\
Email: dengcs@whu.edu.cn
\and \textbf{Wei Liu} \\ Department of Mathematics, Shanghai Normal University, Shanghai 200234, China\\
Email: weiliu@shnu.edu.cn
\and \textbf{Erkan Nane}\\ Department of  Mathematics and Statistics, Auburn University, Alabama 36849, USA \\
Email: ezn0001@auburn.edu
}
\date{}

\begin{document}
\maketitle

\begin{abstract}
The finite time blowup in the almost sure sense of a class of space-time fractional stochastic partial differential equations is discussed. Both the cases of white noise and colored noise are considered. The sufficient and necessary condition between the blowup and Osgood condition is obtained when the spatial domain is bounded. And the sufficient condition for the blowup is obtained when the spatial domain is the whole space. The results in this paper could be regarded as extensions to some results in {\it Foondun and Nualart, 2021}.

\end{abstract}

\maketitle

\section{Introduction}
In this paper, we continue our  research on the blowup of solutions to stochastic partial differential equations (SPDEs). In our previous work \cite{DLN2022}, the finite time blowup {\bf in $L^2$ sense} of solutions to the white or colored noise driven SPDEs  with Bernstein functions of the Laplacian were investigated. Inspired by Foondun and Nualart \cite{foondun-nualart-2019}, we are going to study the finite time blowup of space-time fractional SPDEs {\bf in the almost sure sense} in this paper. Our results in this work could be regarded as extension to some of the theorems in \cite{foondun-nualart-2019}. It should also be mentioned that some blowup results about a class of space-time fractional SPDEs were obtained in Desalegn et al. \cite{DAM2021}, where the fractional operation on the time variable was only considered on the $u_t(x)$ term on the left side of the equations and only the white noise was discussed.
\par
In this paper, we investigate the more general space-time fractional SPDEs and discuss the case of the white noise as well as the colored noise.
\par
We start by considering the space-time fractional SPDEs in a  ball  $B:=B(0,1)$
\begin{equation}\label{eq:dir-additional-source}
\begin{split}
 \partial^\beta_tu_t(x)&=
 -\nu(-\Delta)^{\alpha/2} u_t(x)+I^{1-\beta}[ b(u_t(x))+\sigma\stackrel{\cdot}{F}(t,x)],\quad t>0\quad\text{and}\quad x\in B;\\
 u_t(x)&=0,\quad  x\in \rd\setminus B, \ t>0;\\
 u_t(x)|_{t=0}&=u_0(x).
 \end{split}\end{equation}
Here $\sigma $ is a positive constant,  $-(-\Delta)^{\alpha/2}$ denotes the generator of  $\alpha$-stable  L\'evy process killed upon exiting the ball $B$, $\partial_t^{\beta}$ is the Caputo derivative,  and $I_t^{1-\beta}$ is the fractional integral operator (see Section 2 for the precise definition). The noise $\dot F$, when not space-time white noise is taken to be spatially colored which
is white in time and has a spatial correlation given by the Riesz kernel. That is,
$$E( \stackrel{\cdot}{F}(t,x) \stackrel{\cdot}{F}(s,y)) = \delta_0(t-s)f(x-y),$$
where $f(x) = |x|^{-\eta}, 0 < \eta  < d$. Suppose that
$b $ satisfies the Osgood  condition: for some $a>0$
\begin{equation}\label{osgood-condition}
\int_{a}^\infty\frac{\d s}{b(s)}<\infty.
\end{equation}

When $\eta <\min(2, \beta^{-1})\alpha$, the mild solution of equation \eqref{eq:dir-additional-source} is given in the sense of Walsh \cite{walsh} as follows (see \cite{Mijena-nane-2015} for a motivation to study equations of this  type and \cite{Foondun-mijena-nane-2016} for the proof of the existence of solutions ):

\begin{equation}\label{eq:dir-mild-additional-drift}
\begin{split}
u_t(x)&=\int_{B}G_B({t}, x,y)u_0(y)\,\d y
+ \int_{B}\int_0^t b(u_s(y))G_B({t-s}, x,y)\,\d s\,\d y\\
&\quad+\sigma \int_{B}\int_0^t G_B({t-s},x,y)\,F(\d s\,\d y),
\end{split}
\end{equation}
where $G_B({t}, x,y)$ is the fundamental solution of equation \eqref{eq:dir-additional-source} when $b=0$ and $\sigma=0$ .

The eigenfunctions $\{\phi_n: n\in \mathbb{N}\}$ of fractional Laplacian $-(-\Delta)^{\alpha/2}$ in $B $  form  an orthonormal basis for $L^2(B)$.
Let $E_\beta(z)=\sum_{n=0}^\infty\frac{z^n}{\Gamma(1+n\beta  )}$ denote the Mittag-Leffler function for $\beta\in (0,1)$.
We have an eigenfunction expansion of the kernel
\begin{equation}\label{dir-kernel}
G_B({t}, x,y)=\sum_{n=1}E_\beta(-\mu_nt^\beta)\phi_n(x)\phi_n(y).
\end{equation}
See, for example, Chen et al.   \cite{cmn-12} and  Meerschaert et al. \cite{mnv-09}.

 We will make the following assumption throughout the paper.
\begin{assumption}\label{b-assume-initial-assume}  The function $b : \R \to \R_+$ is nonegative, locally Lipschitz and nondecreasing
on $(0, \infty)$ and the initial condition $u_0(x)$ is nonnegative and continuous.
\end{assumption}
\par
Now, we are ready to present our first theorem and its  proof will be given  in Section 3.
\begin{theorem}\label{thm:domain}
 Suppose that  Assumption \ref{b-assume-initial-assume} holds.
If $\beta>1/2$, and $b$ is convex
and satisfies the Osgood condition \eqref{osgood-condition}, then the solution to \eqref{eq:dir-additional-source} blows up in finite time with positive probability. Conversely, if the solution blows up in finite time with positive probability, then $b$ satisfies the Osgood condition \eqref{osgood-condition}.
\end{theorem}
\begin{remark}
This result could be regarded as an extension of Theorem 1.5 in \cite{foondun-nualart-2019}. Theorem \ref{thm:domain} still holds when the noise is replaced by space-time white one. The proof has similar steps as in the proof of Theorem \ref{thm:domain}, see Bonder and Groisman \cite{bonder-groisman}.
\end{remark}
\par
Now, we turn to the case of the whole space.
Consider the equation with the space-time white noise:
\begin{equation}\label{tfspde}
\begin{split}
 \partial^\beta_tu_t(x)&=-\nu(-\Delta)^{\alpha/2} u_t(x)+I^{1-\beta}_t[b(u_t(x))+\sigma\stackrel{\cdot}{W}(t,x)],\ \ \nu>0, t> 0,\, x\in\R^d;\\
 u_t(x)|_{t=0}&=u_0(x),
 \end{split}
 \end{equation}
where the initial datum $u_0$ is a nonnegative, continuous, and bounded function,
 $-(-\Delta)^{\alpha/2} $ is the fractional Laplacian with $\alpha\in (0,2]$,  $\stackrel{\cdot}{W}(t,x)$ is a  space-time white noise with $x\in \R^d$, and $\sigma$ is a positive constant.

 When $d<\min(2,\beta^{-1})\alpha$, the mild solution of equation \eqref{tfspde} is given in the sense of Walsh \cite{walsh} as follows (see \cite{Mijena-nane-2015}  ]
 $$
        u_t(x)=\int_{\R^d} G(t,x,y)u_0(y)\,\d y
        +\int_0^t\int_{\R^d} G(t-s,x,y)b(u_s(y))\,\d y\,\d s
        +\int_0^t\int_{\R^d} G(t-s,x,y)\,W(\d y\,\d s),
    $$
where $G(t,x,y)$ is the heat kernel of \eqref{tfspde} when $b=0$ and $\sigma=0$.
\par
Now, we are ready to present our second theorem and its  proof  is given  in Section 4.
\begin{theorem}\label{onedim}
  Suppose that Assumption \ref{b-assume-initial-assume} holds. If $b$ satisfies the Osgood  condition, then almost surely, there is no global solution to equation to \eqref{tfspde}.
\end{theorem}

\begin{remark}
This theorem could be regarded as an extension of Theorem 1.4 in \cite{foondun-nualart-2019}. One of the interesting questions is that if the blowup of the solution in finite time in the almost sure sense can indicate that $b$ satisfies the Osgood condition \eqref{osgood-condition}. Unluckily, we have  not found a way to answer it. Meanwhile, we notice that this is also an open question even when the equation is the classical stochastic heat equation, as stated in \cite{foondun-nualart-2019}.
\end{remark}
\par
We notice that there are very some interesting discussions, at the end of Section 1 in \cite{foondun-nualart-2019}, about the non-existence of the solutions to some SPDEs for any $t > 0$, that is the non-existence of the local solution. We would like to make the same claim as Foondun and Nualart did that such non-existence results of the local solutions are not our interests in this paper and our main scope is the non-existence of the global solution. We refer the readers to \cite{foondun-nualart-2019} for more discussions.
\par
The organization of the paper is as follows. Some preliminary results that are needed for the proofs of two theorems are given Section 2. Proofs of Theorem \ref{thm:domain} and Theorem \ref{onedim} are provided  in Sections 3 and 4, respectively.
\section{Preliminaries}

We give some preliminaries that will be essential to the proofs of the main results. Let $\beta\in (0,1)$, $\partial^\beta_t$ is the Caputo fractional derivative which first appeared in \cite{Caputo} and is defined by
$$
\partial^\beta_t u_t(x)=\frac{1}{\Gamma(1-\beta)}\int_0^t \partial_r
u_r(x)\frac{\d r}{(t-r)^\beta} .
$$
For $\gamma>0$, the fractional order integral is defined by
$$
I^{\gamma}_tf(t):=\frac{1}{\Gamma(\gamma)} \int _0^t(t-\tau)^{\gamma-1}f(\tau)\,\d\tau.
$$

We have some useful estimates for the Mittag-Leffler function next.
\begin{lemma}\label{bound}
   We have
   $$
        \inf_{s>0}
        \eup^{\lambda_1s}E_\beta(-\lambda_1s^\beta)
        \,\in\,(0,1].
   $$
\end{lemma}

\begin{proof}
    Recall that (cf.\ \cite[Theorem 4]{Sim14})
    $$
        E_\beta(-x)\geq\frac{1}{1+\Gamma(1-\beta)x},\quad x>0.
    $$
    We have
    $$
        E_\beta(-\lambda_1s^\beta)
        \geq\frac{1}{1+\Gamma(1-\beta)\lambda_1s^\beta}
        \geq C_{\beta,\lambda_1}\,\eup^{-\lambda_1s},
        \quad s>0,
    $$
    where
    $$
        C_{\beta,\lambda_1}:=\inf_{s>0}
        \frac{\eup^{\lambda_1s}}
        {1+\Gamma(1-\beta)\lambda_1s^\beta}\,\in\,(0,1].
    $$
    This implies the desired assertion.
\end{proof}

Denote by $p^\alpha(t,x,y)$ the
 transition density of a symmetric
stable process $X_t$ in $\rd$ of index $\alpha\in (0, 2]$. Denote by $E_t$
the inverse $\beta$-stable subordinator.
Then the time-changed
Brownian motion $X_{E_t}$ has a transition density
$$
    G(t,x,y)=G(t,x-y)=\E p^\alpha(E_t, x,y)=\int p^\alpha(s,x,y)\,\P(E_t\in\d s).
$$

\begin{lemma}\label{low--gen-d}
Let $d\geq 1$. We have
    $$
        \inf_{x\in B(0,1),r\in(0,1)}
        \int_{B(0,1)}G(r,x,y)\,\d y>0
    $$
\end{lemma}

\begin{proof}
    If $x\in B(0,1)$ and $r\in(0,1)$, it is easy
    to verify that
    $$
        B\left((1-r^{\beta/\alpha}/2)x,
        r^{\beta/\alpha}/2\right)\,\subset\,
        B(0,1)\cap B(x,r^{\beta/\alpha}).
    $$
    By \cite[Lemma 2.1 (a)]{foondun-nane-2017} there is some $c_0>0$ such that
    $$
       G(t,x,y)\geq c_0t^{-\beta d /\alpha}\quad
        \text{whenever $|x-y|\leq t^{\beta/\alpha}$}.
    $$
    Then for any $x\in B(0,1)$ and $r\in(0,1)$,
    \begin{align*}
        \int_{B(0,1)}G(r,x,y)\,\d y
        &\geq\int_{B(0,1)\cap B(x,r^{\beta/\alpha})}
        c_0r^{-\beta d/\alpha}\,\d y\\
        &\geq c_0r^{-\beta d/\alpha}
        \int_{B\left((1-r^{\beta/\alpha}/2)x,
        r^{\beta/\alpha}/2\right)}\d y\\
        &=c_0r^{-\beta d/\alpha}
        \frac{\pi^{d/2}}{\Gamma(d/2+1)}
        \left(\frac{r^{\beta/\alpha}}{2}\right)^d\\
        &=\frac{c_0\pi^{d/2}}{2^d\Gamma(d/2+1)}.\qedhere
    \end{align*}
\end{proof}

\begin{lemma}[Lemma 1 in \cite{Mijena-nane-2015}]\label{Lem:Green1} For $d < 2\a,$
$$
\int_{{\R^d}}[G(t,x)]^2\,\d x  =C^\ast t^{-\beta d/\a},
$$
where $C^\ast>0$ is a constant depending only on 
$d,\alpha,\beta$.
\end{lemma}

In the following, we consider
 $$
 g_{\beta}(t,x):=\int_0^t\int_{\R^d}G(t-r,x-y)\,W(\d r\,\d y).
$$
This is the random part of the mild solution to \eqref{tfspde} when $\sigma\equiv 1$.

\begin{lemma}\label{lem-LIL}Suppose $d<\min(2,\beta^{-1})\a$.
For each fixed $x\in \rd$, almost surely,
$$
\limsup_{t\to\infty}\frac{g_{\beta}(t,x)}{\sqrt{Kt^{1-\beta d/\alpha}\log \log t}}=1.
$$
\end{lemma}

\begin{proof}
Since $\{g_{\beta}(t,x), \ t\ge 0\}$ is a Gaussian process for each fixed $x$,
we need to calculate the variance of $g_{\beta}(t,x)$. By the Plancharel theorem
 and using Lemma \ref{Lem:Green1},
$$
\E[g_{\beta}(t,x) g_{\beta}(t,x)]=\int_{\rd}\int_0^{t}[G(t-r,x-y)]^2
\,\d r\, \d y=\frac{C^\ast}{1-\beta d/\alpha}\, t^{1-\beta d/\alpha}.
$$
Hence by equation (5) in Lai  \cite{Lai-1974} we get the result using the fact that the variance of $ g_{\beta}(t,x)$ is given by Lemma \ref{Lem:Green1}
$$
\limsup_{t\to\infty}\frac{g_{\beta}(t,x)}{\sqrt{Kt^{1-\beta d/\alpha}\log \log t}}=1,
$$
where $K:=\frac{2C^\ast}{1-\beta d/\alpha}$.
\end{proof}

\begin{proposition}[Proposition 2 \cite{Mijena-nane-2015}]\label{Prop:MomentEst1}
Suppose $d<\min(2,\beta^{-1})\a$, and $k\ge 2$. Then there
exists $c>1$ such that the following moment estimates for time increments and spatial increments hold.

\noindent (i). For $s\le t$, 
$$
c^{-1}|t-s|^{\left(1-\frac{\beta d}{\a}\right)\frac k2}\le \E\left[|g_{\beta}(t,x)-g_{\beta}(s,x)|^k\right]\le c|t-s|^{\left(1-\frac{\beta d}{\a}\right)\frac k2}.
$$
(ii). For  $x,\, y\in{\R^d}$,
$$
c^{-1}|x-y|^k\le\E\left[|g_{\beta}(t,x)-g_{\beta}(t,y)|^k\right]\le c|x-y|^{\min\left\{\left(\frac{\a-\beta d}{\beta}\right)^{-},2\right\}\frac k2}.
$$
\end{proposition}

 As a consequence of Proposition  \ref{Prop:MomentEst1} and classical Garsia's lemma we have the following estimate.
 \begin{proposition}\label{prop-est-in-n}For all $k\ge 2$, there exists $A_k>0$ such that for any $n\ge 1$,
 $$
 \E\bigg[ \sup_{s,t\in [n, n+2], x,y \in B(0, 1)}|g_{\beta}(t,x)-g_{\beta}(s,y)|^k\bigg]\leq A_k 2^{k(1-\beta d/\alpha)/2}.
 $$
 \end{proposition}

Let $\Psi(t):=\sqrt{Kt^{1-\beta d/\alpha}\log \log t}$.
Now using this lemma, we get the next proposition.

\begin{proposition}\label{prop-lim-increment-n}
Almost surely,
$$
\sup_{s,t\in [n, n+2], x,y \in B(0,1)}\frac{|g_{\beta}(t,x)-g_{\beta}(s,y)|}{\Psi(n)}\to 0,\ \ \mathrm{as}\ \ n\to\infty.
$$

\end{proposition}

\begin{proof} Proposition \ref{prop-est-in-n} implies that for  $k\ge 2 \alpha/(\alpha-\beta d)$

$$
\E \bigg[ \sum_{n=1}^\infty  \sup_{s,t\in [n, n+2], x,y \in B(0,1)}\frac{|g_{\beta}(t,x)-g_{\beta}(s,y)|^k}{\Psi(n)^k}\bigg]\leq \sum_{n=1}^\infty \frac{A_k 2^{k(1-\beta d/\alpha)/2}}{\Psi(n)^k}<\infty.
$$
This gives the desired result.

\end{proof}

 Using Proposition \ref{prop-lim-increment-n} we get the following result.

 \begin{proposition}\label{prop-infty-sequence}
 Almost surely, there exists a sequence $t_n\to\infty$ such that

 $$
 \inf_{h\in[0,1], x\in B(0,1)} g_{\beta}(t_n +h,x)\to \infty \ \ \mathrm{as}\ \ n\to\infty.
 $$

 \end{proposition}

 \begin{proof}
  Fix $x_0\in B(0,1)$. Choose $\omega$ such that Lemma \ref{lem-LIL} and Proposition \ref{prop-lim-increment-n} hold. Then we get
  \begin{align*}
  \inf_{h\in[0,1], x\in B(0,1)} g_{\beta}(t +h,x)& =
    g_{\beta}(t ,x_0)+ \inf_{h\in[0,1], x\in B(0,1)} \big[g_{\beta}(t +h,x)-g_{\beta}(t ,x_0)\big]\\
  &\ge g_{\beta}(t ,x_0)+ \inf_{h\in[0,1], x\in B(0,1)} -\big[|g_{\beta}(t +h,x)-g_{\beta}(t ,x_0)|\big]\\
  &= \frac{g_{\beta}(t ,x_0)}{\Psi(t)}\Psi(t)+ -\sup_{h\in[0,1], x\in B(0,1)} \frac{\big[|g_{\beta}(t +h,x)-g_{\beta}(t ,x_0)|\big]}{\Psi([t])}\Psi([t])\\
  &=\min \{\Psi(t), \Psi([t])\}\bigg[\frac{g_{\beta}(t ,x_0)}{\Psi(t)}-\sup_{h\in[0,1], x\in B(0,1)} \frac{\big[|g_{\beta}(t +h,x)-g_{\beta}(t ,x_0)|\big]}{\Psi([t])} \bigg].
 \end{align*}
Now using  Lemma \ref{lem-LIL} and Proposition \ref{prop-lim-increment-n}, we can choose a suitable 
sequence  $t_n $  to finish the proof.
 \end{proof}

\section{Proof of Theorem \ref{thm:domain}}

\begin{proof}[Proof of Theorem \ref{thm:domain}]
(1) To prove the first assertion, we shall borrow  an idea used in the proof of Theorem 1.3 in  Foondun and Nualart \cite{foondun-nualart-2019}.

\smallskip\noindent
a) Set
$$Y_t:=\int_{B}u(t,x)\phi_1(x)\,\d x.$$
From  equation \eqref{dir-kernel},  we can easily get
$$
\int_B  G_B({t}, x,y)\phi_1(x)\,\d x=  E_\beta(-\mu_1t^\beta)  \phi_1(y).
$$
It is a well-know fact that $\phi_1(x)>0$ for $x\in B$.

Then after using the stochastic Fubini theoem and the decomposition above for the heat kernel $G_B(t-s,x,y)$, from the mild solution \eqref{eq:dir-mild-additional-drift}, we get
$$
Y_t=E_\beta(-\lambda_1 t^{\beta})Y_0
+\int_0^t  E_\beta(-\lambda_1 (t-s)^{\beta})\int_B b(u_s(y))\phi_1(y)\,\d y\,\d s
+ \int_0^t E_\beta(-\lambda_1 (t-s)^{\beta})\int_B\phi_1(y)
\,F(\d y\,\d s).
$$
Since $b$ is convex, using the Jensen inequality, we can get
$$
\int_B b(u_s(y))\phi_1(y)\,\d y \geq  b(Y_s)\geq
b\big(Y_s-E_\beta(-\lambda_1 s^{\beta})Y_0\big),
$$
where the second inequality follows from the assumptions
that $b$ is nondecreasing and $Y_0\geq0$. We also have
 $$
 \int_0^t \int_B\phi_1(y)\,F(\d y\,\d s)=\sqrt{\kappa}B_t,
 $$
where $B_t$ is a Brownian motion and
$$\kappa:=\int_{B\times B}\phi_1(y)\phi_1(z)
|y-z|^{-\eta}\,\d y\,\d z,$$
where $\eta$ is the Riesz kernel exponent.

Therefore we obtain
$$
Y_t\geq E_\beta(-\lambda_1 t^{\beta})Y_0
+\int_0^t  E_\beta(-\lambda_1 (t-s)^{\beta})b\big(Y_s-E_\beta(-\lambda_1 s^{\beta})Y_0\big)
\,\d s
+ \sqrt{\kappa}\int_0^t E_\beta(-\lambda_1 (t-s)^{\beta})\,
\d B_s.
$$
By Lemma \ref{bound}, there exists $c=c(\beta,\lambda_1)>0$
such that
$$
Y_t\geq E_\beta(-\lambda_1 t^{\beta})Y_0
+c\,\int_0^t  \eup^{-\lambda_1(t-s)}
b\big(Y_s-E_\beta(-\lambda_1 s^{\beta})Y_0\big)
\,\d s
+ \sqrt{\kappa}\int_0^t E_\beta(-\lambda_1 (t-s)^{\beta})\,
\d B_s.
$$

\smallskip\noindent
b) To use the comparison principle, consider
$$
Z_t=E_\beta(-\lambda_1 t^{\beta})
\,Z_0
+c\int_0^t \eup^{-\lambda_1(t-s)} b(Z_s-
E_\beta(-\lambda_1 s^{\beta})Z_0
)\,\d s
+ \sqrt{\kappa}
\,\int_0^tE_\beta(-\lambda_1 (t-s)^{\beta})\,\d B_s
$$
with $Z_0=Y_0$. We have
\begin{align*}
\eup^{\lambda_1t}Z_t&=\eup^{\lambda_1t}E_\beta(-\lambda_1 t^{\beta})
\,Z_0+c\int_0^t \eup^{\lambda_1s} b(Z_s-
E_\beta(-\lambda_1 s^{\beta})Z_0
)\,\d s\\
&\quad+ \sqrt{\kappa}\,\eup^{\lambda_1t}
\int_0^tE_\beta(-\lambda_1 (t-s)^{\beta})\,\d B_s.
\end{align*}
This implies that
\begin{align*}
\d Z_t&=\d\{\eup^{-\lambda_1t}(\eup^{\lambda_1t}Z_t)\}
=-\lambda_1Z_t\,\d t+
\eup^{-\lambda_1t}\,\d\{\eup^{\lambda_1t}Z_t\}\\
&= -\lambda_1Z_t\,\d t+c\,b(Z_t-E_\beta(-\lambda_1 t^{\beta})Z_0)\,\d t
+\sqrt{\kappa}\,\d B_t\\
&\quad+\left(\lambda_1E_\beta(-\lambda_1 t^{\beta})Z_0
+\frac{\d E_\beta(-\lambda_1 t^{\beta})}{\d t}\,
Z_0
\right)\,\d t\\
&\quad+\sqrt{\kappa}\left(
\int_0^t
\left[\lambda_1E_\beta(-\lambda_1 (t-s)^{\beta})
+\frac{\d E_\beta(-\lambda_1 (t-s)^{\beta})}{\d t}
\right]\,\d B_s\right)\,\d t.
\end{align*}
Setting $U_t:=Z_t-E_\beta(-\lambda_1 t^{\beta})Z_0$,
we get
\begin{equation}\label{ueq}
\d U_t= -\lambda_1U_t\,\d t+c\,b(U_t)\,\d t
+\sqrt{\kappa}\,\d B_t+\sqrt{\kappa} \xi(t)\,\d t,
\end{equation}
where
$$
\xi(t):=\int_0^t
\left[\lambda_1E_\beta(-\lambda_1 (t-s)^{\beta})
+\frac{\d E_\beta(-\lambda_1 (t-s)^{\beta})}{\d t}
\right]\,\d B_s.
$$

\smallskip\noindent
c) Consider
$$
    \d V_t= -\lambda_1V_t\,\d t+c\,b(V_t)\,\d t
+\sqrt{\kappa}\,\d B_t
$$
with $V_0=U_0$. By the Feller test, $V_t$ explodes
in finite time a.s. This means that there exists
(deterministic) $T<\infty$ such that
\begin{equation}\label{posi}
    \P\left(\text{there exists $S=S(\omega)\leq T$ such that}
    \;\; \lim_{t\uparrow S}V_t=\infty\right)>0.
\end{equation}

\smallskip\noindent
d) Let
$$
    h(t):=\lambda_1E_\beta(-\lambda_1 t^{\beta})
    +\frac{\d E_\beta(-\lambda_1 t^{\beta})}{\d t}
    =\lambda_1E_\beta(-\lambda_1 t^{\beta})-\beta\lambda_1
    t^{\beta-1}E_\beta'(-\lambda_1 t^{\beta}),\quad t>0,
$$
and $h(0):=0$. Since $\beta>1/2$, $h\in L^2_{\text{loc}}(\R_+)$. Note that
$$
    \xi(t)=\int_0^th(t-s)\,\d B_s.
$$
By Lemma \ref{unif} in Appendix,
$$
    \tilde{B}_t:=B_t+\int_0^t\xi(s)\,\d s,\quad 0\leq t\leq T,
$$
is a Brownian motion under the weighted probability
measure $R_T\P$, where
$$
    R_T:=\exp\left[
    -\int_0^T\xi(s)\,\d B_s-\frac12\int_0^T|\xi(s)|^2\,\d s
    \right].
$$
Rewrite \eqref{ueq} as
$$
    \d U_t= -\lambda_1U_t\,\d t+c\,b(U_t)\,\d t
+\sqrt{\kappa}\,\d \tilde{B}_t.
$$
Then we know that the distribution of $(U_t)_{0\leq t\leq T}$
under $R_T\P$ coincides with that of $(V_t)_{0\leq t\leq T}$
under $\P$. By \eqref{posi}, we conclude that
$$
    (R_T\P)\left(\text{there exists $S=S(\omega)\leq T$
    such that}
    \;\; \lim_{t\uparrow S}U_t=\infty\right)>0.
$$
This implies
\begin{equation}\label{blow}
    \P\left(\text{there exists $S=S(\omega)\leq T$
    such that}
    \;\;  \lim_{t\uparrow S}U_t=\infty\right)>0.
\end{equation}
By the comparison principle,
$Y_t\geq Z_t=U_t+E_\beta(-\lambda_1 t^{\beta})Z_0$.
Consequently, \eqref{blow} holds with
$U_t$ replaced by $Y_t$, and the proof is now finished.

\bigskip

(2) Next, we follow the line of  the proof of Theorem 1.3 in Foondun and Nualart \cite{foondun-nualart-2019} (with crucial changes) to prove the second assertion. Let
  $$
      T:=\sup\left\{t\geq0\,:\,\sup_{x\in B(0,1)}|u_t(x)|
      <\infty\right\}.
  $$
  Since the solution blows up in finite time with positive probability,
we can find a set $A$
 satisfying $P(A) > 0$ such that for any $\omega \in A$, we have $T (\omega) < \infty$. In the following
  we write $T=T(\omega)<\infty$ (drop the variable $\omega$).
  Note that $T$ is the blowup time.

 Recall that ($\sigma= 1$)
  \begin{equation*}
\begin{split}
u_t(x)&=
\int_{B(0,1)}G_B(t,x,y)u_0(y)\,\d y+ \int_0^t\int_{B(0,1)} b(u_s(y))G_B({t-s}, x,y)\,\d y\,\d s\\
&\quad+\int_0^t\int_{B(0,1)} G_B({t-s},x,y)\,F(\d y\,\d s)\\
&=:\sum_{i=1}^3I_i(t,x).
\end{split}
\end{equation*}
Since the initial value $u_0$ is bounded, we find that
$$
  |I_1(t,x)|\leq \int_{B(0,1)}G_B(t,x,y)|u_0(y)|\,\d y
  \leq \|u_0\|_\infty\int_{B(0,1)}G_B(t,x,y)\,\d y
  =\|u_0\|_\infty,
 \quad  \forall t\in[0,T],\,\forall x\in B(0,1).
$$
Set $Y_t:=\sup_{x\in B(0,1)}u_t(x)$. Since $b$
is nondecreasing, one has
$$
    I_2(t,x)\leq \int_0^t\int_{B(0,1)} b(Y_s)G_B({t-s}, x,y)\,\d y\d s
    =\int_0^tb(Y_s)\,\d s,\quad \forall t\in[0,T).
$$
Noting that $I_3$ is continuous almost surely, there is
some $M>0$ such that
$$
    \sup_{t\in[0,T],x\in B(0,1)}|I_3(t,x)|\leq M.
$$
Substituting these estimates into the first formula, we get
$$
    u_t(x)\leq \|u_0\|_\infty+\int_0^tb(Y_s)\,\d s
    +M,\quad \forall t\in[0,T),\,\forall x\in B(0,1).
$$
Taking supremum over $x\in B(0,1)$,
$$
    Y_t\leq \widetilde{M}+\int_0^tb(Y_s)\,\d s,\quad \forall
    t\in[0,T),
$$
where $\widetilde{M}:=\|u_0\|_\infty+M$. Consider
 $$
    Z_t= \widetilde{M}+\int_0^tb(Z_s)\,\d s
$$
with $Z_0=Y_0$. By the comparison principle,
$Y_t\leq Z_t$. Since $Y_t$ blows
up at time $t=T$, so does $Z_t$. This immediately
implies that $b$ satisfies the so-called
Osgood condition by the classical 
ODE theory (cf. \cite{Osg98}).
\end{proof}

\section{Proof of theorem \ref{onedim}}

\begin{proof}[Proof of Theorem \ref{onedim}]
    Recall the mild formulation
    $$
        u_t(x)=\int_{\R^d} G(t,x,y)u_0(y)\,\d y
        +\int_0^t\int_{\R^d} G(t-s,x,y)b(u_s(y))\,\d y\,\d s
        +\int_0^t\int_{\R^d} G(t-s,x,y)\,W(\d y\,\d s).
    $$
    Let $\{t_n\}$ be a sequence  such that $t_n\rightarrow\infty$ as $ n\rightarrow\infty$  and
    Proposition  \ref{prop-infty-sequence}
    holds. Since $b$ and $u_0$ are nonnegative,
    \begin{align*}
        u_{t+t_n}(x)
        &=\int_{\R^d} G(t+t_n,x,y)u_0(y)\,\d y
        +\int_0^{t+t_n}\int_{\R^d} G(t-s,x,y)
        b(u_s(y))\,\d y\,\d s\\
        &\quad+\int_0^{t+t_n}\int_{\R^d}G(t-s,x,y)
        \,W(\d y\,\d s)\\
        &\geq
        \int_{\R^d} G(t+t_n,x,y)u_0(y)\,\d y+\int_0^t\int_{\R^d} G(t-s,x,y)b(u_{s+t_n}(y))\,
        \d y\,\d s\\
        &\quad+\int_0^{t+t_n}\int_{\R^d} G(t+t_n-s,x,y)\,
        W(\d y\,\d s)\\
        &=:\sum_{i=1}^3I_i(t,x,n).
    \end{align*}
    Let
    $$
        Y_{t,n}:=\inf_{y\in B(0,1)}u_{t+t_n}(y).
    $$
    By Lemma \ref{low--gen-d}  for any $x\in B(0,1)$,


    \begin{align*}
        I_2(t,x,n)&\geq\int_0^t\int_{B(0,1)} G(t-s,x,y)b(u_{s+t_n}(y))\,
        \d y\,\d s\\
        &\geq \int_0^t b(Y_{s,n})
        \int_{B(0,1)}G(t-s,x,y)\,
        \d y\,\d s\\
        &\geq c\int_0^t b(Y_{s,n})\,\d s.
    \end{align*}
    Note that for any $x\in B(0,1)$ and $t\in(0,1)$
    $$
        I_3(t,x,n)=g_\beta(t_n+t,x)\geq
        \inf_{h\in[0,1],x\in B(0,1)}g_\beta(t_n+h,x).
    $$
    Then we conclude that for
    any $x\in B(0,1)$, $t\in(0,1)$ and $n\in\nat$
    $$
        u_{t+t_n}(x)\geq
        \inf_{h\in[0,1],x\in B(0,1)}g_\beta(t_n+h,x)
        +c\int_0^t b(Y_{s,n})\,\d s
    $$
    Taking infimum over $x\in B(0,1)$, we know
    that for any $t\in(0,1)$ and $n\in\nat$,
    $$
        Y_{t,n}\geq
        \inf_{h\in[0,1],x\in B(0,1)}g_\beta(t_n+h,x)
        +c\int_0^t b(Y_{s,n})\,\d s
    $$
    Consider for $t\in(0,1)$ and $n\in\nat$
    \begin{equation}\label{zblowup}
        Z_{t,n}=
        \inf_{h\in[0,1],x\in B(0,1)}\bigg[ \int_{\R^d} G(t+t_n,x,y)u_0(y)\,\d y +g_\beta(t_n+h,x)\bigg]
        +c\int_0^t b(Z_{s,n})\,\d s
    \end{equation}
    By the comparison principle,
    $Y_{t,n}\geq Z_{t,n}$.

\bigskip

Suppose that the solution does not blow up in finite
time. Then
for any $n\in\nat$, $Y_{t,n}$ does
not blow up. This means that the blow up
time of $Z_{t,n}$ has to be greater than $1$.
By the classical ODE theory (cf. \cite{Osg98}) 
and \eqref{zblowup},
$$
    \int_{\inf_{h\in[0,1],x\in B(0,1)}g_\beta(t_n+h,x)}
    ^\infty \frac{\d s}{b(s)}>1.
$$
Letting $n\rightarrow\infty$, it follows from
the Osgood condition that the left hand side
tends to zero. This leads to a contradiction
and thus finishes the proof.
\end{proof}

\section*{Appendix}

Let $B_t$ be a standard Brownian motion on $\R$.
Let $h\in L^2_{\text{loc}}(\R_+)$, and set for $t\geq0$
$$
    \xi(t):=\int_0^th(t-s)\,\d B_s,\quad
    R_t:=\exp\left[-\int_0^t\xi(s)\,\d B_s
    -\frac12\int_0^t|\xi(s)|^2\,\d s
    \right].
$$

\begin{lemma}\label{unif}
    For any $T>0$, $\{R_t\}_{0\leq t\leq T}$ is a uniformly
    integrable martingale.
\end{lemma}

\begin{proof}
    For $n\in\N$, let
    $$
        \tau_n:=\inf\left\{t\geq0\,;\,\int_0^t|\xi(s)|^2
        \,\d s\geq n\right\}.
    $$
    By the Girsanov theorem, for any $t\in(0,T]$,
    $$
        \tilde{B}_s:=B_s+\int_0^s\xi(s)\,\d s,
        \quad 0\leq s\leq t,
    $$
    is a Brownian motion under the
    weighted probability measure $R_{t\wedge\tau_n}\P$.
    Using the elementary inequality that
    $|x-y|^2\leq2x^2+2y^2$ for $x,y\in\R$, we obtain that
    for $s\in[0,t]$
    \begin{align*}
        \E_{R_{t\wedge\tau_n}\P}\left[|\xi(s\wedge\tau_n)|^2
        \right]
        &=\E_{R_{t\wedge\tau_n}\P}\left[\left|
        \int_0^{s\wedge\tau_n}h(s\wedge\tau_n-r)\,\d \tilde{B}_r
        -\int_0^{s\wedge\tau_n}h(s\wedge\tau_n-r)\xi(r)\,\d r
        \right|^2\right]\\
        &\leq2\,\E_{R_{t\wedge\tau_n}\P}\left[\left|
        \int_0^{s\wedge\tau_n}h(s\wedge\tau_n-r)\,\d \tilde{B}_r
        \right|^2\right]
        +2\,\E_{R_{t\wedge\tau_n}\P}\left[\left|
        \int_0^{s\wedge\tau_n}h(s\wedge\tau_n-r)\xi(r)\,\d r
        \right|^2\right]\\
        &=:2I_1+2I_2.
    \end{align*}
    By the It\^{o} isometry, one has
    \begin{align*}
        I_1&=\E_{R_{t\wedge\tau_n}\P}\left[
        \int_0^{s\wedge\tau_n}|h(s\wedge\tau_n-r)|^2\,\d r
        \right]
        =\E_{R_{t\wedge\tau_n}\P}
        \left[\int_0^{s\wedge\tau_n}|h(r)|^2\,\d r\right]\\
        &\leq\int_0^{s}|h(r)|^2\,\d r
        \leq\int_0^{t}|h(r)|^2\,\d r=:A_t.
    \end{align*}
    It follows from the H\"{o}lder inequality that
    \begin{align*}
        I_2&\leq\E_{R_{t\wedge\tau_n}\P}\left[
        \int_0^{s\wedge\tau_n}|h(s\wedge\tau_n-r)|^2
        \,\d r\,
        \cdot\,
        \int_0^{s\wedge\tau_n}|\xi(r)|^2\,\d r\right]\\
        &=\E_{R_{t\wedge\tau_n}\P}\left[
        \int_0^{s\wedge\tau_n}|h(r)|^2
        \,\d r\,
        \cdot\,
        \int_0^s|\xi(r\wedge\tau_n)|^2\,\d r\right]\\
        &\leq\int_0^s|h(r)|^2
        \,\d r\,
        \cdot\,\E_{R_{t\wedge\tau_n}\P}\left[
        \int_0^s|\xi(r\wedge\tau_n)|^2\,\d r\right]\\
        &\leq A_t\,\cdot\,
        \int_0^{s}\E_{R_{t\wedge\tau_n}\P}
        \left[|\xi(r\wedge\tau_n)|^2\right]
        \,\d r.
    \end{align*}
    Then we get
    $$
        \E_{R_{t\wedge\tau_n}\P}\left[
        |\xi(s\wedge\tau_n)|^2\right]
        \leq2A_t+2A_t\,\cdot\,
        \int_0^{s}\E_{R_{t\wedge\tau_n}\P}
        \left[|\xi(r\wedge\tau_n)|^2\right]
        \,\d r,\quad 0\leq s\leq t,
    $$
    which, together with the Gronwall inequality, implies that
    $$
        \E_{R_{t\wedge\tau_n}\P}\left[|\xi(s\wedge\tau_n)|^2
        \right]
        \leq2A_t\,\eup^{2A_t s}.
    $$
    This yields that for $t\in(0,T]$
    \begin{align*}
        \E\left[R_{t\wedge\tau_n}\log R_{t\wedge\tau_n}
        \right]&=\E_{R_{t\wedge\tau_n}\P}\left[
        -\int_0^{t\wedge\tau_n}\xi(s)\,\d \tilde{B}_s
        +\frac12\int_0^{t\wedge\tau_n}|\xi(s)|^2\,\d s
        \right]\\
        &=\E_{R_{t\wedge\tau_n}\P}
        \left[\frac12\int_0^{t\wedge\tau_n}|\xi(s)|^2\,\d s\right]\\
        &=\frac12\int_0^{t}
        \E_{R_{t\wedge\tau_n}\P}|\xi(s\wedge\tau_n)|^2\,\d s\\
        &\leq A_t\int_0^{t}
        \eup^{2A_ts}\,\d s
        =\frac12\big(\eup^{2A_tt}-1\big).
    \end{align*}
    Since $\lim_{n\rightarrow\infty}\tau_n=\infty$,
    it holds from the Fatou lemma that
    $$
        \E\left[R_t\log R_t\right]\leq
        \frac12\big(\eup^{2A_tt}-1\big).
    $$
    Thus,
    $$
        \sup_{t\in[0,T]}\E\left[R_t\log R_t\right]\leq
        \frac12\big(\eup^{2A_TT}-1\big)<\infty,
    $$
    and this completes the proof.
\end{proof}

\end{document}